\documentclass[11pt]{amsart}
\usepackage{amsmath} 
\usepackage{amssymb} 
\usepackage{amsfonts} 
\usepackage{graphicx}
 
 \newtheorem{thm}{Theorem}[section] 
 \newtheorem{cor}[thm]{Corollary} 
 \newtheorem{lem}[thm]{Lemma} 
 \newtheorem{prop}[thm]{Proposition} 
 \theoremstyle{definition}

 \theoremstyle{remark} 
 \newtheorem{rem}[thm]{Remark}

\numberwithin{equation}{section} 
\numberwithin{figure}{section} 
 
\newcommand{\CC}{\mathbb{C}}
\newcommand{\DD}{\mathbb{D}} 
\newcommand{\TT}{\mathbb{T}} 
\newcommand{\EE}{\mathbb{E}} 
\newcommand{\PP}{\mathbb{P}}

\newcommand{\fr}{\frac}
\newcommand{\ee}{\epsilon}

\newcommand{\iy}{\infty}
\newcommand{\vp}{\varphi}
\newcommand{\diam}{\operatorname{diam}}

\begin{document} 
\bibliographystyle{alpha}

\title[L\'evy-Loewner hulls and random growth] 
{Rescaled L\'evy-Loewner hulls and random growth} 
\author[Johansson] 
{Fredrik Johansson} 
 
\address{Johansson: Department of Mathematics\\ 
The Royal Institute of Technology\\  
S -- 100 44 Stockholm\\ 
SWEDEN} 
 
\email{frejo@math.kth.se} 
 
\thanks{Research supported by grant KAW 2005.0098 from 
the Knut and Alice Wallenberg Foundation.} 
 
\author[Sola] 
{Alan Sola} 
 
\address{Sola: Department of Mathematics\\ 
The Royal Institute of Technology\\ 
S -- 100 44 Stockholm\\ 
SWEDEN} 
\email{alansola@math.kth.se}

 
  
\begin{abstract} 
We consider radial Loewner evolution driven by unimodular L\'evy processes. 
We rescale the hulls of the evolution 
by capacity, and prove that the weak limit of the rescaled hulls exists. 
We then study a random growth model obtained 
by driving the Loewner equation with a compound Poisson process.
The process involves two real parameters: the intensity of the underlying 
Poisson process and a localization parameter of the Poisson kernel which determines 
the jumps. A particular choice of parameters yields a growth process similar to 
the Hastings-Levitov $\rm{HL}(0)$ model. We describe the asymptotic behavior 
of the hulls with respect to the parameters, showing that growth tends to 
become localized as the jump parameter increases. We obtain deterministic 
evolutions in one limiting case, and Loewner evolution driven by a unimodular 
Cauchy process in another. We show that the Hausdorff dimension of the 
limiting rescaled hulls is equal to $1$. Using a different type of compound 
Poisson process, where the Poisson kernel is replaced by 
the heat kernel, as driving function, we recover one case of the 
aforementioned model and $\rm{SLE}(\kappa)$ as limits. 
\end{abstract} 
 
\maketitle 
 
\section{Introduction and main results} 
\subsection{Introduction}
The {\it Loewner differential equation} was first derived by Loewner  
in 1923 in a paper on coefficient problems for univalent functions and  
has proved to be one of the most powerful tools in the theory of univalent  
functions. This differential equation parametrizes functions mapping  
conformally onto various reference domains (usually the unit disk or 
the upper half-plane) minus a single slit  
in terms of a continuous unimodular function. On the other hand, given a  
unimodular function (which will be called {\it driving  
function} in this context), one obtains from the Loewner equation a univalent  
function, which needs not be a single-slit mapping.  

Loewner's equation has also proved to be of use in mathematical physics.
For instance, in the paper \cite {CM}, Carleson and Makarov use the Loewner 
equation to study a model related to Diffusion Limited Aggregation (DLA).
In the last few years, solutions to the Loewner equation corresponding to  
random driving functions have attracted a great amount of interest in 
mathematics as well as physics. For instance, if standard Brownian motion on 
the unit circle is chosen as driving function, one obtains 
{\it Schramm-Loewner evolution} or SLE for short. SLE arises naturally 
in connection with scaling limits of various discrete models from statistical
physics and provides a tool to treat many aspects of that field rigorously. 
We refer the reader to the book \cite{Lawler-book} and the paper 
\cite{Rohde_Schramm} for basic results on SLE. Recently, the Loewner 
equation together with other, more general, stochastic processes has
attracted some interest, see \cite{Rohde_Stable} and the references therein. 
\subsection{Main results} 
Our paper is organised as follows. Beginning in Section \ref{LP} we consider 
the Loewner equation driven by general L\'evy processes. This yields a family 
of conformal maps mapping the exterior disk onto the complements of a family 
of growing 
compact connected sets, the so-called hulls of the evolution. We take on 
the view of the Loewner equation as a mapping from the Skorokhod space of 
right-continuous functions with left limits, and establish the continuity of 
this mapping in Proposition \ref{cont_of_LM}. We rescale the compacts by capacity and by using reversibility of L\'evy processes, we prove a form of convergence of the rescaled hulls in Theorem \ref{existence}. 

In the final section, 
we study 
evolutions corresponding to a particular choice of L\'evy process of pure 
jump type involving two real parameters $r$ and $\lambda$. We 
consider random variables $X^r$ with densities proportional to harmonic 
measure at different points whose distance to the growing cluster 
is controlled by $r$. The variables are then added 
at the occurrence times of a Poisson process with intensity $\lambda$.
For the particular choice $r=0$ (corresponding to harmonic measure at the 
point at infinity), the 
evolution we obtain resembles the Hastings-Levitov model $\rm{HL}(0)$ in Laplacian growth and may 
be viewed as a random perturbation of that model. We study our models' dependence on the 
parameters, 
and show that the hulls we obtain become more regular with increasing $r$ (which 
corresponds to moving the base points of harmonic measure closer to the hull) in the 
sense that the 
growth becomes more localized and the hulls tend to grow fewer large branches. In the 
limit, with 
$\lambda$ fixed, we obtain deterministic hulls, see Proposition \ref{det-hulls}. Another passage to the limit 
(which involves both 
parameters) yields Loewner evolution driven by the Cauchy process, see Proposition \ref{cauchy}. The hulls generated 
by this process still exhibit branching behavior. Replacing the Poisson kernel as
the density of the variables $X^r$ by the heat kernel of the circle, which is more 
localized, we obtain Schramm-Loewner evolution $\rm{SLE}(\kappa)$ in a particular 
limiting case that involves both the intensity of the Poisson process and the 
localization parameter of the heat kernel. At the other extreme, when the heat kernel 
becomes spread out on the circle, we recover the perturbed Hastings-Levitov 
model. 
Finally, in Theorem \ref{Hdim}, we show that the Hausdorff dimension of the 
scaling limit of the hulls (in the sense of Theorem \ref{existence}) is $1$ 
for all $r$ and large enough $\lambda$, as is the case in the $\rm{HL}(0)$ 
model. 
\section{Preliminaries}
\subsection{The Loewner equation and classes of conformal mappings}
Let $\DD$ denote the unit disk and let $\xi$ be a
function of the real variable $t$ taking values on the unit circle.
The \emph{radial} version of Loewner's equation is 
\begin{equation} 
\partial_t g_t(z)=g_t(z)\frac{\xi(t)+g_t(z)}{\xi(t)-g_t(z)},
\quad \textrm{with}\quad g_0(z)=z,
\label{Loewnereq} 
\end{equation} 
for $0\leq t<\infty$.
For each $z \in \DD$, the equation (\ref{Loewnereq}) is 
well-defined up to 
\begin{equation*}\tau_z = \sup\{t:\inf_{0 \le s \le t}|g_s(z)-\xi(s)| > 0\},\end{equation*} 
that is, the first time the denominator in the right-hand side
of (\ref{Loewnereq}) is zero. We let 
\begin{equation}
K_t=\{z: \tau_z \le t, \, z \in \DD\} 
\end{equation}
denote the growing compact \emph{hulls}. 
It is well-known that there exists a unique solution $g_t$ to 
(\ref{Loewnereq}) if $\xi$ is integrable (see \cite[Section 3.4]{Duren} for 
a general discussion). 
For each $t$, the function $g_t$ is a conformal map
from $\DD \setminus K_t$ onto $\DD$. 

We note that $g_t$ fixes the origin and 
that $g'_t(0)=e^t$. Often, it is the inverse of $g_t$, that is, 
the conformal map of $\DD$ onto $\DD \setminus K_t$, that we would like to 
study. This function may be obtained as the solution of a partial differential equation 
(with given initial value) involving the driving function, namely,
\begin{equation}
\partial_tf_t(z)=z\partial_z f _t(z) \frac{z+\xi(t)}{z-\xi(t)}, \quad f_0(z)=z.
\label{partialloewner1}
\end{equation}

Next, let
\begin{equation*} 
\Delta=\{z \in \mathbb{C}_{\infty}:|z|>1\}, 
\end{equation*}
be the exterior of the unit disk. 
It is easily verified that the map $g^*_t=1/g_t(1/z)$ also satisfies 
(\ref{Loewnereq}), this time with driving function $\overline{\xi}$. 
The solution $g^*_t$ maps $\Delta \setminus K_t^*$ conformally onto 
$\Delta$, where 
\begin{equation*}K_t^*=\{z: 1/z \in K_t\}.\end{equation*} 
We note that $g^*_t$ has an expansion of the form
\begin{equation*}g^*_t(z)=e^{-t}z+\cdots\end{equation*}
at infinity. Again, we obtain the inverse mappings by solving 
the partial differential equation (\ref{partialloewner1}). Hence, the equations (\ref{Loewnereq}) 
and (\ref{partialloewner1}) can be 
used to describe growing hulls both in $\DD$ and $\Delta$. In this paper, we
consider evolutions in the exterior disk.

One may also consider the equation (\ref{Loewnereq}) in the exterior disk 
for negative $t$. This leads to the notion of \emph{whole-plane} Loewner 
evolution. More precisely, for a driving function $\vp$ 
defined for $-\iy<t<\iy$, we consider solutions of 
\begin{equation}
\label{Whole-plane}
\partial_t g_t(z) = g_t(z)\fr{\vp(t)+g_t(z)}
{\vp(t)-g_t(z)}, 
\end{equation}
with initial condition
\begin{equation}
\quad \lim_{t \to -\iy}e^{t}g_t(z)=z.  
\end{equation}
The associated hulls, defined in the same manner as before, grow from the 
origin and become arbitrarily large as $t \to \iy$. As in the previous cases, 
the inverse mappings $f_t$ satisfy the related partial differential equation
\begin{equation}
\partial_tf_t(z)=z\partial_z f_t(z) \frac{z+\vp(t)}
{z-\vp(t)}, 
\label{partialloewner}
\end{equation}
with initial condition
\begin{equation}
\lim_{t \rightarrow -\iy}e^{-t}f_t(z)=z.
\label{partialinit}
\end{equation}
We shall occasionally refer to these different Loewner equations by using the 
prefixes $\DD, \, \Delta$, and $\CC$.  

We have seen that the Loewner equation generates conformal mappings, and 
we now set down some notation for the different classes of mappings that 
arise. The class $\mathcal{S}$ consists of univalent functions 
$f:\DD\rightarrow \CC$ with 
power series expansions of the form
\begin{equation*}f(z)=z+\sum_{n=2}^{\infty}a_n z^n.\end{equation*}
These functions map the unit disk onto simply connected domains 
containing the origin. The Loewner equation was in fact first introduced in 
connection with coefficient problems for $\mathcal{S}$. Next, we let $\Sigma$ 
denote the class of univalent 
functions $f$ in $\Delta$ with expansions at infinity of the form 
\begin{equation*} 
f(z)=z+b_0+\sum_{n=1}^{\infty}\frac{b_n}{z^n}. 
\end{equation*} 
A function in $\Sigma$ maps the exterior disk onto the complement 
of a compact simply connected set. Finally, the subclass $\Sigma'$ consists of those 
$f \in \Sigma$ whose omitted set contains the origin. We note that 
$\mathcal{S}$ and $\Sigma'$ are related via the inversion mapping discussed 
earlier. We endow the classes $\mathcal{S}$ and $\Sigma$ with the 
topology induced by uniform convergence on compact subsets. We note 
that $\Sigma'$ is compact with respect to this topology. We refer the reader 
to \cite{Duren} for a thorough treatment of univalent functions.

\subsection{L\'evy processes and Skorokhod space}
A (real-valued) stochastic process $X=(X(t))_{t \in [0,\iy)}$ defined on some 
probability space $(\Omega,\mathcal{F}, \PP)$ is said to 
be a {\it L\'evy process} if $X(0)=0$ almost surely, if the process has 
independent and stationary increments, and if, for all $\ee>0$, it holds that 
\begin{equation*}\lim_{t \rightarrow 0}\PP(|X(t)|>\ee)=0.\end{equation*}
Clearly, standard Brownian motion is an example of a L\'evy process. The
Poisson process and compound Poisson processes provide us with additional 
examples of L\'evy processes. It is well-known that the sample paths of 
standard Brownian motion are continuous almost surely. The sample paths of 
general L\'evy processes do not have this property. It is true, however, 
that the sample paths of a L\'evy process are right-continuous with left 
limits (RCLL) (more precisely, any L\'evy process admits a modification with 
these properties). 

There is a natural topology, known as the {\it Skorokhod topology}, on 
spaces of RCLL functions on the real line. Let us first consider the 
case of Skorokhod spaces on an interval $[0,T]$. The class 
$\Lambda=\Lambda[0,T]$ consists of continuous non-decreasing functions on 
$[0,T]$ with $\lambda(0)=0$ and $\lambda(T)=T$.  
We set 
\begin{equation}
\|\lambda\|_{\Lambda}=\sup_{s\neq t}\left|\log\left(
\frac{\lambda(t)-\lambda(s)}{t-s}\right)\right|.
\end{equation}
The Skorokhod metric on $D[0,T]$, the space of real-valued RCLL functions on 
$[0,T]$, is then given by
\begin{equation}
d(\vp,\psi)=\inf\{\epsilon>0: \exists \lambda \,\, \textrm{such that} \,\,
\|\vp-\psi\circ \lambda\|_{\iy}<\epsilon, 
\|\lambda\|_{\Lambda}<\epsilon\}.
\end{equation} 
The space $D[0,T]$ is a complete and separable metric space when equipped 
with this metric. Note that the restriction of $d$ to 
the subspace of continuous functions induces the usual uniform 
topology on $C[0,T]$. The following lemma (see \cite{Billingsley}) shows that
the functions in $D[0,T]$ are reasonably well-behaved:
\begin{lem}
\label{uniformity}
For each $\epsilon>0$ and $\vp \in D[0,T]$ there exist a finite number 
$N=N(\epsilon)$ of points $0=t_0<t_1\cdots<t_N=T$ such that
\begin{equation*}
\sup\{|\vp(s)-\vp(t)|:s,t\in I_j\}<\epsilon
\end{equation*}
for  $j=0,1,\ldots N-1$ and $I_j=[t_{j},t_{j+1})$.
\end{lem}
We see that Skorokhod functions can have at most a countable 
number of discontinuities. The above lemma, together with the 
fact that $\|\lambda\|_{\Lambda}<\epsilon$ implies $\|\lambda-t\|_{\infty} 
<2T\epsilon$, can be used to show that convergence in $D[0,T]$ in the 
Skorokhod metric implies convergence in $L^1[0,T]$. In fact, using 
that a function in $D[0,T]$ can be approximated uniformly by 
a finite step function, one can prove the 
following stronger result (we omit the details).
\begin{lem}
\label{Skorokhod-appr}
For each $\epsilon>0$ and each $\vp \in D[0,T]$, there exists a
$\delta>0$ such that $d(\vp,\psi)<\delta$ implies that for all 
$t \in [0,T]$,
\begin{equation*}\int_0^t|\vp(s)-\psi(s)|ds<\epsilon t.\end{equation*}
\end{lem}

We often wish to consider the space $D[0,\iy)$ of 
RCLL functions on the real line. The Skorokhod topology is defined as follows 
in this case (see \cite{Jacod_book}). We say that a sequence 
$\{\vp_i\}_{i=1}^{\iy}$ of functions in $D[0,\iy)$ converges to 
$\vp \in D[0,\iy)$ if there is a sequence of continuous, strictly 
increasing functions $\{\lambda_i\}_{i=1}^{\iy}$ such that $\lambda_i(0)=0$, 
$\lambda_i(t) \rightarrow \iy$ as $t \rightarrow \iy$, and
\begin{equation}
\left\{ \begin{array}{c}
\sup_{t \in [0,\iy)} |\lambda_i(t)-t|\rightarrow 0\\
\sup_{t\leq N} |\vp_i(\lambda_i(t))-\vp(t)|\rightarrow 0, \quad
\forall N \in \mathbb{N}.\end{array}\right.
\end{equation}
The Skorokhod topology on $D(-\iy,0]$ is defined in an analogous manner.

We will need the following lemma which is contained in \cite[Corollary VII.3.6]{Jacod_book}.
\begin{lem}\label{levy-conv}
Let $(X_n(t))_{t \in [0,\iy)}$ and $(X(t))_{t \in [0,\iy)}$ be L\'evy processes. Then, as $n \to \iy$, the law of $X_n$ converges weakly to that of $X$ in the space of probability measures on $D[0, \iy)$ if and only if the random variables $X_n(1)$
converge in distribution to $X(1)$.
\end{lem}

We refer the reader to \cite{Bertoin} for more background information on 
L\'evy processes. The Skorokhod space is treated in depth in 
\cite{Jacod_book}.

Throughout this paper, we shall always assume that probability measures are 
defined on the appropriate Borel sigma algebra. For a metric space $S$, we 
let $\Pi(S)$ denote the space of Borel probability measures on $S$ with the 
topology induced by weak convergence.

\section{Rescaled Loewner evolution driven by L\'evy processes}\label{LP}
This section is dedicated to the study of Loewner evolution driven 
by general L\'evy processes, and we establish several facts that 
we shall use later on in a more specific setting. The main result of 
this section is the weak convergence of Loewner evolution 
driven by L\'evy processes after rescaling by capacity. 
\subsection{Continuity of the Loewner mapping}
Let $D_{\mathbb{T}}[0,\iy)$ denote the space of unimodular RCLL functions 
on $[0,\iy)$, equipped with the Skorokhod topology. The space
$L^1_{\mathbb{T}}[0,T]$ consists of all measurable unimodular 
functions on $[0,T]$. We note that any function in $D_{\mathbb{T}}[0,\iy)$
belongs to $L^1_{\mathbb{T}}[0,T]$ for any choice of $T$.

Consider the time $t=T$ solution $f_T$ to the $\Delta$-Loewner equation 
(\ref{partialloewner}) with driving function $\vp \in L^1_{\TT}[0,\iy)$ 
and initial condition $f_0(z)=z$.
We define a mapping 
\begin{equation}
\mathcal{L}_T: L^1_{\mathbb{T}}[0,T] \rightarrow \Sigma'
\end{equation}
by setting
\begin{equation}
\mathcal{L}_T[\vp](z)=e^{-T}f_T(z).
\end{equation}
Usually, we shall consider the mapping $\mathcal{L}_T$ restricted to the 
Skorokhod space $D_{\mathbb{T}}[0,\iy)$. We next show that the mappings 
$\mathcal{L}_T$ are well-behaved. 
\begin{prop}\label{cont_of_LM}
For each $T>0$, the Loewner mapping 
\begin{equation*}\mathcal{L}_T:L^1_{\mathbb{T}}[0,T]
\rightarrow \Sigma'\end{equation*} 
is continuous.
\end{prop}
We will follow the approach of Bauer. He obtains similar
results for the chordal Loewner equation viewed as a mapping from $C[0, \iy)$, 
with the uniform topology, to a certain space of 
conformal mappings of the upper half-plane (see \cite{Bauer} for details; see
also \cite[Section 3.4]{Duren}). 
\begin{proof}[Proof of Proposition \ref{cont_of_LM}]
It is well-known (see for instance \cite[Section 4.2]{Lawler-book}) that the conformal 
mapping $f_T$ can be obtained by considering the backward flow
\begin{equation*}
\partial_t h_t(z)=-h_t(z)\frac{\vp(T-t)+h_t(z)}{\vp(T-t)-h_t(z)},
\quad h_0(z)=z
\end{equation*}
for $0\leq t\leq T$, and we have $\mathcal{L}_T[\vp](z)=e^{-T}h_T(z)$.
In what follows, we set
\begin{equation*}
v(t,\vp,z)=-z\frac{\vp(T-t)+z}{\vp(T-t)-z}.
\end{equation*}
With this notation, the equation determining the backward flow 
becomes
\begin{equation}
\partial_t h_t(z)=v(t,\vp,h_t(z)),\quad h_0(z)=z.
\label{backflow}
\end{equation}

Let $\epsilon>0$ be given, and suppose $\vp, \psi \in L^1_{\mathbb{T}}[0,T]$.
We fix $z \in \Delta$ with $z=re^{i\theta}$, and we let 
$u_1=u_1(t)$ denote the solution to (\ref{backflow}) with $\vp$ as 
driving function. Similarily, we let $u_2$ denote the solution corresponding 
to $\psi$. Moreover, we write $\dot{u}_i$ for the derivative of $u_i$
with respect to $t$. We now 
show that $u_2$ is an approximate solution to $\dot{u}_1(t)=v(t,\vp,u_1)$.

We compute that 
\begin{equation*}
\dot{u}_2(t)-v(t,\vp,u_2)=
-2u_2^2(t)\frac{\vp(T-t)-\psi(T-t)}{(\vp(T-t)-u_2(t))(\psi(T-t)-u_2(t))}.
\end{equation*}
An integration then yields
\begin{multline*}
\left|u_2(t)-u_2(0)-\int_0^t v(s,\vp,u_2)ds\right|\\
\leq 2 \int_0^t|u_2(s)|^2\frac{|\vp(T-s)-\psi(T-s)|}
{|\vp(T-s)-u_2(s)||\psi(T-s)-u_2(s)|}ds.
\end{multline*}
Standard growth estimates show that $|u_2(s)|\leq 4e^Tr$. 
Next, we calculate that
\begin{equation*}\frac{d}{ds}\textrm{Re}[\log u_2(s)]>0,\end{equation*}
which means that points in $\Delta$ move away from the unit circle 
under the backward flow. This implies that
\begin{equation*}|\vp(T-s)-u_2(s)|\geq |u_2(0)|-1=r-1,\end{equation*} 
and a similar estimate holds 
for the other factor. Using these estimates, we obtain
\begin{multline*}
\left|u_2(t)-u_2(0)-\int_0^t v(s,\vp,u_2)ds\right|\\
\leq
32e^T\frac{r^2}{(r-1)^2}\int_0^t|\vp(T-s)-\psi(T-s)|ds.
\end{multline*}
The integral on the right-hand side can be made arbitrarily small
if the functions $\vp$ and $\psi$ are chosen so that 
$\|\vp-\psi\|_1\leq 1/32\, e^{-T}\epsilon$. We then get
\begin{equation*}
\left|u_2(t)-u_2(0)-\int_0^t v(s,\vp,u_2)ds\right|\leq
\frac{r^2}{(r-1)^2}\epsilon.
\end{equation*}
Noting that $u_1(0)=u_2(0)=z$, and that $u_1$ is an exact solution 
to the problem (\ref{backflow}), we see that
\begin{equation*}
|u_1(t)-u_2(t)|\leq \frac{r^2}{(r-1)^2}\epsilon+
\int_0^t|v(s,\vp,u_1)-v(s,\vp,u_2)|ds.
\end{equation*}
We wish to estimate the integrand in terms of the difference $|u_1-u_2|$. 
A computation shows that
\begin{equation*}
  |\partial_z v(s,\vp,z)|\leq \frac{(r+1)^2}{(r-1)^2},
\end{equation*}
and we obtain that
\begin{equation*}
  |u_1(t)-u_2(t)|\leq \frac{r^2}{(r-1)^2}\epsilon
  +\frac{(4e^Tr+1)^2}{(r-1)^2}\int_0^t|u_1(s)-u_2(s)|ds
\end{equation*}
Finally, an application of Gr\"onwall's lemma shows that
\begin{equation*}
|u_1(t)-u_2(t)|\leq 
\frac{r^2}{(r-1)^2}\exp\left\{\frac{(4e^Tr+1)^2}{(r-1)^2}T\right\}
\cdot \epsilon.
\end{equation*}
We see that this estimate is uniform for $z=r e^{i \theta}$ confined 
to an arbitrary compact set in $\Delta$. Since
$|u_1(T)-u_2(T)|=|\mathcal{L}_T[\varphi](z)-\mathcal{L}_T[\psi](z)|$, 
the proof is now complete.
\end{proof}
We recall that convergence in the Skorokhod topology on $[0,\iy)$
implies convergence in $L^1{[0,T]}$ and obtain the following corollary.
\begin{cor}
\label{continuity}
For each $T>0$, the Loewner mapping 
\begin{equation*}\mathcal{L}_T:D_{\mathbb{T}}[0,\iy) 
\rightarrow \Sigma'\end{equation*} 
is continuous.
\end{cor}

\subsection{Weak convergence of rescaled hulls}
Let $X=(X(t))_{t \in [0,\iy)}$ be a L\'evy process
and let $U$ be a uniformly distributed random variable on the unit circle, 
independent of $X$.
We introduce a driving function for the $\Delta$-Loewner equation by 
setting
\begin{equation}
  \xi(t)=Y(t)=\exp\{i(U+X(t))\}, \quad t \in [0,\infty).
    \label{levydriving}
\end{equation}
We then let $f_t$, $t>0$, 
be the solution to (\ref{partialloewner}) with driving function $\xi(t)$.

In this section we study rescalings of the hulls $K_t=\mathbb{C}\setminus f_t(\Delta)$. Recall that the capacity of the hulls $K_t$ is given by $\textrm{cap}(K_t)=e^t$.
In view of our application below is natural to normalize the hulls so that they have capacity $1$ for 
each $t$; this also means the diameter of the hulls stays bounded (both from above and away from $0$). 
We therefore introduce the rescaled hulls
\begin{equation}
\tilde{K}_t=\frac{1}{\textrm{cap}(K_t)}K_t=e^{-t}K_t.
\label{rescaledhull}
\end{equation}
On the level of mappings, this normalization places the conformal 
mapping $\tilde{f}_t$ of $\Delta$ onto $\mathbb{C} \setminus \tilde{K}_t$ in the 
class $\Sigma'$. Let $\PP_t$ denote the law of $\tilde{f}_t$. 
\begin{thm}
We have 
\begin{equation*}
\PP_t \rightarrow \PP_{\infty},
\end{equation*}
as $t \to \iy$ in the sense of weak convergence of probability measures on $\Sigma'$.
\label{existence}
\end{thm}
The limiting measure $\PP_{\iy}$ will correspond to a whole-plane Loewner evolution evaluated at $t=0$.

Our approach requires that we introduce a suitable driving function for negative $t$. This involves time-reversal of the driving L\'evy process.

In general, L\'evy processes are reversible in the following sense. Suppose
$V$ is a L\'evy process on $[0, T]$. We define the time-reversed process 
$(\check{V}(t))_{t \in [0,T]}$ by setting 
\begin{equation}
\check{V}(t)=\left\{ \begin{array}{cc}0 & t=0\\
                           V((T-t)-) -V(T-)& 0<t<T\\
			   -V(T-) & t=T \end{array}. \right.
\end{equation}
Here we use the notation $V(s-)$ for the left-limit of $V$ at the
point $s$. The process $(\check{V}(t))_{t \in [0,T]}$ is again a L\'evy process and has the distribution of $(-V(t))_{t \in [0,T]}$ (see \cite{Protter} for details). 

Since we work with L\'evy processes on the unbounded interval $[0,\iy)$, we will need to slightly modify the 
construction above in order to obtain a suitable driving process 
for negative $t$: Let $X$ be a given L\'evy process on $[0, \iy)$ with law 
$\mu:=\mu_X \in \Pi(D[0, \iy))$. We define a reflected and time reversed process
$\check{X}=(\check{X}(t))_{t \in  (-\iy, 0]}$, by setting 
\begin{equation}
\check{X}(t)=\left\{ \begin{array}{cc}0 & t=0\\
                  		  -X((-t)-) & t<0 \end{array}. \right.
\end{equation}
We view $\check{X}$ as a random element of the space $D(-\iy, 0]$, 
with law $\check{\mu}:=\mu_{\check{X}} \in \Pi(D(-\iy, 0])$. 
\begin{lem}
\label{time-reversal} Fix $T > 0$ and let $X_1$ and $X_2$ be independent 
L\'evy processes with laws $\mu \in \Pi(D[0, \iy))$ and
$\check{\mu} \in \Pi(D(-\iy, 0])$ respectively. Let $U_1$ and $U_2$ be two random variables uniformly distributed on $[-\pi, \pi]$, chosen independently of each other, as well as $X_1$ and $X_2$.
Then the two processes
\begin{equation*}
Y_1(t):=\exp\{i(U_1+X_1(t)\}, \quad 0 \le t \le T
\label{posproc}
\end{equation*}
and
\begin{equation*}
Y_2(t):=\exp\{i(U_2+X_2(t-T))\}, \quad 0 \le t \le T
\label{negproc}
\end{equation*}
have the same law as elements of $D_{\TT}[0,T]$.
\end{lem}
\begin{proof}
By reversibility, on $[0,T)$, it holds that
\begin{eqnarray*}
X_1(\cdot) &\stackrel{d}{=}& -(X_1((T-\cdot)-) - X_1(T-)),
\end{eqnarray*}
and hence on $[0,T)$ 
\begin{eqnarray*}
Y_1(\cdot) &\stackrel{d}{=}&
\exp\{i[U_1-\{X_1((T-\cdot)-) - X_1(T-)\}]\} \\
&\stackrel{d}{=}&\exp\{i(U_2-X_1((T-\cdot)-))\}.
\end{eqnarray*}
However, on $[0,T)$ we have
\begin{eqnarray*}
\exp\{i(U_2-X_1((T-\cdot)-))\} &\stackrel{d}{=}& 
\exp\{i(U_2+X_2(\cdot-T))\} \\
 &=&Y_2(\cdot),
\end{eqnarray*}
and since $X_2(0)=0$ almost surely, this implies that the processes are 
equally distributed on the whole interval.
\end{proof}
We keep the notation of the previous lemma and let $\nu \in \Pi(D_{\TT}[0, \iy))$ and $\check{\nu} \in \Pi(D_{\TT}(-\iy, 0])$ denote the laws of $Y_1$ and 
$Y_2$ respectively. We let $K_t$ denote the hulls corresponding to $Y_1$.

Recall that $\mathcal{L}_t$ is continuous for each $t \ge 0$, and that $\PP_t$ is the law of the unique map $f_t \in \Sigma'$ such that $e^{-t}K_t=\CC \setminus f_t(\Delta)$. It follows that the family
of probability measures $(\PP_t)_{t \in[0,\iy)}$ on the class $\Sigma'$ can be written
\begin{equation*}
\PP_t=\nu \circ \mathcal{L}_{t}^{-1}.
\end{equation*}
We define a mapping $L_{\iy}: D_{\TT}(-\iy,0] \rightarrow \Sigma'$ by evaluating the solution to the 
whole-plane Loewner equation at $t=0$. An important step in the proof 
of Theorem \ref{existence} below is to show that $L_{\iy}$ is continuous on 
$D_{\TT}(-\iy,0]$, so that the composition
\begin{equation*}\PP_{\iy}=\check{\nu}\circ L_{\iy}^{-1}\end{equation*} 
is another probability measure on $\Sigma'$.
\begin{proof}[Proof of Theorem \ref{existence}]
We shall use the notation  
\begin{equation*}
w(t,\vp, z; f):=z \partial_z f_t(z)\frac{z+\vp(t)}{z-\vp(t)}
\end{equation*}
for the right hand side of the Loewner partial differential equation.
Let $\{t_n\}_{n=0}^{\iy}$ be an increasing sequence such that $t_0=0$ and 
$t_n \to \iy$ as $n \to \iy$. We write $\PP_{n}:=\PP_{t_n}$.

For an arbitrary $\vp \in D_{\TT}(-\iy, 0]$, we let $f^n_t$ be the
(unique) solution to $\partial_t f^n_t = w(t,\vp, z; f^n)$ for 
$-t_n < t \le 0$ with the initial value $f^n_{-t_n}(z)=e^{-t_n}z$,
and we set
\begin{equation*}
L_n[\vp]:= f^n_0.
\end{equation*}
It is easy to see, using the same arguments as in the previous section, that 
the mappings $L_n$ are continuous on $D_{\TT}(-\iy, 0]$, so that $\check{\nu} \circ L_n^{-1}$ is a probability measure. In view of Lemma 
\ref{time-reversal},
we then have $\check{\nu} \circ L_n^{-1}
=\nu \circ \mathcal{L}_{t_n}^{-1}$ so that we may write 
$\PP_{n}=\check{\nu}\circ L_n^{-1}$. Hence, 
the rescaled hull $\tilde{K}_{n}:=e^{-t_n}K_{t_n}$ has the distribution of 
$\mathbb{C}\setminus f(\Delta)$, where $f$ is chosen according to 
$\check{\nu} \circ L_n^{-1}$. 

We will now show that $\{L_n\}$ converges uniformly on $D_{\TT}(-\iy,0]$. 
Fix some compact $E \subset \Delta$ and take $m > n$ so that $t_m < t_n$.  
Let $\vp \in D_{\TT}(-\iy,0]$ be arbitrary and let $h_{t}^{n}$ solve 
\begin{equation*}\partial_t h_{t}^{n} = w(t, \vp, z;h^{n}),\quad -t_n<t\leq 0,\end{equation*} 
with initial value 
\begin{equation*}h_{-t_n}^n(z)=z.\end{equation*}
Keeping our previous notation, we denote by $f_{t}^{m}$ the solution to 
$\partial_t f_{t}^{m} = w(t, \vp, z;f^{m})$ for $-t_m < t \le 0$, with initial 
value $f_{-t_m}^{m}(z)=e^{-t_m}z$. Using a well-known property of solutions 
to the Loewner equation, we see that
\begin{equation*}L_n[\vp]=e^{-t_n}z \circ h_{0}^{n},\end{equation*}
and
\begin{equation*}L_m[\vp]=f_{-t_n}^{m}\circ h_{0}^{n}.\end{equation*}
We now apply \cite[Lemma 4.20]{Lawler-book} 
in order to conclude that
\begin{equation*}
\sup_{z \in E}|L_n[\vp](z) - L_m[\vp](z)| \le C \cdot e^{-t_n},
\end{equation*}   
for a constant $C < \iy$ independent of $\vp$. 

Thus, the sequence $L_n$ converges 
uniformly on $D_{\TT}(-\iy, 0]$ to the continuous map $L_{\iy}$, given by $L_{\iy}[\vp]=f_0^{\iy}$, where $f_t^{\iy}$ solves 
\begin{equation*}\partial_tf_t^{\iy}=w(t,\vp,z; f^{\iy}), \quad  -\iy < t \le 0\end{equation*} 
with $\lim_{t \to -\iy}e^{-t}f_{t}(z)=z$ as initial condition.

It follows, by bounded convergence for instance, that the measures $\PP_n=\check{\nu}\circ L_n^{-1}$ converge to the measure $\PP_{\iy}=\check{\nu}\circ L_{\iy}^{-1}$ in $\Pi(\Sigma')$. This completes the proof.
\end{proof}
\begin{rem}
In \cite{Rohde_Stable}, Chen and Rohde perform a similar rescaling by capacity in the
 setting of chordal Loewner evolutions in a half-plane driven by $\alpha$-stable
processes. There, it is shown that the limit is trivial in the sense that
$\lim_{t \to \iy}\PP(K_{t^2}/t \cap \{ \textrm{Im}\, z > \ee \} \neq \emptyset)=0$ for all $\ee$. The intuition is that the rescaled hulls are spread out along the real 
line. The results above show that the limiting $\alpha$-stable hulls are 
nontrivial in our setting.
\end{rem}
\section{Loewner evolution driven by a certain compound Poisson process} 
We now consider a random growth model defined using Loewner evolutions corresponding to a 
particular choice of L\'evy driving process. 

This process is of pure jump 
type, and its definition involves two real parameters $r$ and $\lambda$. 
Our motivation for the construction below is to define a somewhat flexible 
growth model similar to the $\rm{HL}(0)$ model by using a L\'evy process.
Let $0 \le r < 1$ be fixed. We introduce a random variable $X^r$ with density
proportional to harmonic measure in the unit disk at the point $z=r$. Explicitly,
the density is given by the Poisson kernel,
\begin{equation*} 
f_{X^r}(\theta)=\frac{1}{2 \pi}\frac{1-r^2}{1-2r \cos\theta+r^2}.
 \label{poiss_densfunct}
\end{equation*}
In the sequel, we shall think of $r$ as a parameter.
Next, we let $N=(N(t))_{t \in [0,\iy)}$ be a Poisson process with intensity 
$\lambda$. We take $\{X_j\}_{1}^{\infty}$ to be a sequence of i.i.d.\
random variables distributed according to the law of $X^r$ and independent of 
$N$. 
We then define the compound Poisson process $Y=(Y(t))_{t\in [0,\iy)}$ by 
setting $Y(0)=0$ and,
\begin{equation}
Y(t)=\sum_{j=1}^{N(t)}X_j,\quad t>0.
\label{procdef}
\end{equation}
Note that $Y(t)=0$ up to the first occurrence time of $N$.
Finally, the driving function $\xi=\xi_{r, \lambda}$ is defined by 
\begin{equation} 
\xi(t)=\exp\{i Y(t)\}, \quad t\in [0, \infty).
\label{xidef} 
\end{equation}  
The process $Y$, and hence the function $\xi$, depend on the two parameters 
$r$ and $\lambda$, but we usually suppress this dependence to keep our
notation as simple as possible.

We shall consider the random growth model given by the $\Delta$-Loewner 
evolution corresponding to $\xi$.
The results on rescaling presented in the previous section 
apply to the thus obtained hulls, providing us with a way to
grow random compact sets with infinitely many branchings. 

We now give some intuition and motivation for the construction of our model.
The $\rm{HL}(0)$ growth model is obtained by composing
random, independent, and uniform rotations of a {\it fixed} conformal map
\begin{equation*}h_{\delta}:\Delta \rightarrow \Delta \setminus [1,1+\delta].\end{equation*}
The general Hastings-Levitov model $\rm{HL}(\alpha)$ and its connections
to physical growth models are discussed in 
\cite{HL}.
One way of thinking about our evolutions is to compose 
random, independent rotations of {\it random} mappings of the type
\begin{equation*}h_{\delta_k}:\Delta \rightarrow \Delta \setminus [1,1+\delta_k].\end{equation*}
At step $k$, a slit of random length $\delta_k$ is mapped to the 
growing hull by the ``previous'' composition of mappings $f_{k-1}$ which is 
obtained as a solution to the Loewner equation, evaluated at 
the $(k-1)^{th}$ occurrence time of the Poisson process $N$. This is 
achieved by precomposing with a rotation of $h_{\delta_k}$.
In other words, if $T_k$ are the occurence times, we have 
\begin{equation*}
f_{k}:=f_{T_k}=h_{\theta_1,\delta_1} \circ \cdots \circ h_{\theta_k,\delta_k},
\end{equation*} 
where $h_{\theta, \delta}(z)=e^{i\theta}h_{\delta}(e^{-i\theta}z)$ and the 
$\theta_j$:s are i.i.d.\ and distributed according to the law of $X^r$.  

The length $\delta_k$ of the slit whose image is added is a function 
of the waiting time $\tau_k$ of the Poisson process underlying the 
evolution. In fact, by considering the explicit expression for 
the basic slit mapping $h_{\delta}$ given in \cite{HL}, one finds that
\begin{equation}
\delta_k=2e^{2\tau_k}(1+\sqrt{1-e^{-2\tau_k}})-2,
\label{deltadist}
\end{equation}
where $\tau_k$ is exponentially distributed with parameter $\lambda$.
A calculation shows that $\delta \to 0$ in probability as $\lambda \to \iy$; we
shall usually think of $\lambda$ as large.
 
The parameter $r$ determines the distribution of the rotations.
When $r=0$, the rotations are uniformly distributed  
on $\TT$, as in the $\rm{HL}(0)$ case. By conformal invariance 
of harmonic measure, this corresponds to adding arcs to the growing hull
according to harmonic measure at infinity. For other values of $r$ the 
$k^{th}$ arc is added according to harmonic measure at the point 
$f_{k-1}(\xi(\tau_{k-1})/r)$, that is, ``close'' to the tip of the 
$(k-1)^{th}$ arc. We shall see below that, in a way, 
the parameter $r$ can be thought of as localizing the growth. As $r$ 
increases to $1$, the jumps of the process $Y$ tend to get smaller, and growth 
of the hulls is then concentrated to fewer ``arms''. In the limit, we obtain 
deterministic hulls. 

\subsection{Asymptotic behaviour of $(r, \lambda) \mapsto K_t^{r, \lambda}$}
In this section, we shall write out both parameters in subscript  
only when we wish to stress depedence on both parameters.
We begin by fixing $\lambda$ and investigate how the choice of $r$ affects 
the process $Y$ by considering the limiting cases. 

Our analysis relies on the Fourier coefficients of the process $\xi$; the 
convergence of Fourier coefficients allows us to deduce the convergence 
in law of the corresponding processes in the Skorokhod space using Lemma \ref{levy-conv} (see also Chapter VII in \cite{Jacod_book} for an in-depth discussion). 
For $t$ fixed, we write $a_{Y(t)}$ for the Fourier 
transform of the law of $Y(t)$; recall that $Y$ is real-valued. 
Since $\xi(t)=\exp(i Y(t) )$, we obtain, by considering the periodization of 
the density of $Y(t)$, that the Fourier coefficients of the law of $\xi$ 
are given by $a_{Y(t)}(n), \, n\in \mathbb{Z}$. Since $Y(t)$ is a sum of 
$N(t)$ i.i.d. variables that are independent of $N(t)$, it follows that
\begin{equation*}
a_{Y(t)}(n)=g_{N(t)}(a_{X^r}(n)).
\end{equation*}
Here, $g_{N(t)}(s)=\exp(\lambda t(s-1))$ is the generating function of $N(t)$, 
and $a_{X^{r}}(n)$ is the $n^{th}$ Fourier coefficient of the law of 
$X^{r}$. Recalling the Fourier representation of 
the Poisson kernel, we see that 
\begin{equation}
a_{X^{r}}(n)=\left\{ \begin{array}{cc}1 & n=0\\ 
r^{|n|} & n \neq 0.\end{array}\right.
\end{equation} 
Hence, we have 
$a_{Y(t)}(0)=1$ and
\begin{equation}
a_{Y(t)}(n)=\exp \left[ \lambda t \left( r^{|n|} - 1 \right) 
\right],\quad n \neq 0.
\label{Ycoeffs}
\end{equation} 
We see that
\begin{equation*}
\lim_{r \to 1}a_{Y(t)}(n)=1,
\end{equation*}
for all $t$ and $n$. This means that the distribution of the random variable 
$\xi(t)$ degenerates to the Dirac measure at the point $1$ in this limiting case. 
Similarily, we see that the limit
\begin{equation*}
\lim_{r \rightarrow 0}a_{Y(t)}(n)=e^{-\lambda t}, \quad n \neq 0,
\end{equation*} 
corresponds to jumps chosen uniformly on $\TT$.
\begin{figure}
\centering
\includegraphics[height=0.24 \textheight, width=0.48 \textwidth]{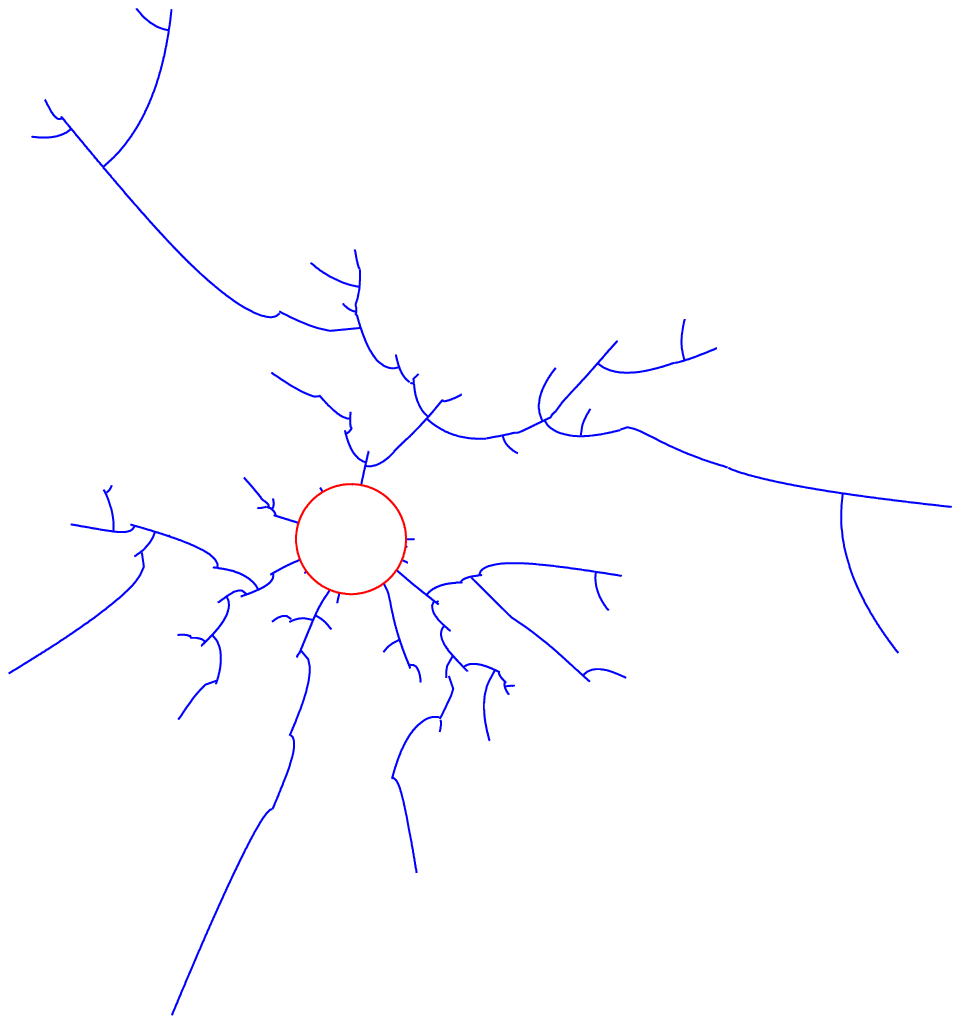}
\includegraphics[height=0.24 \textheight, width=0.48 \textwidth]{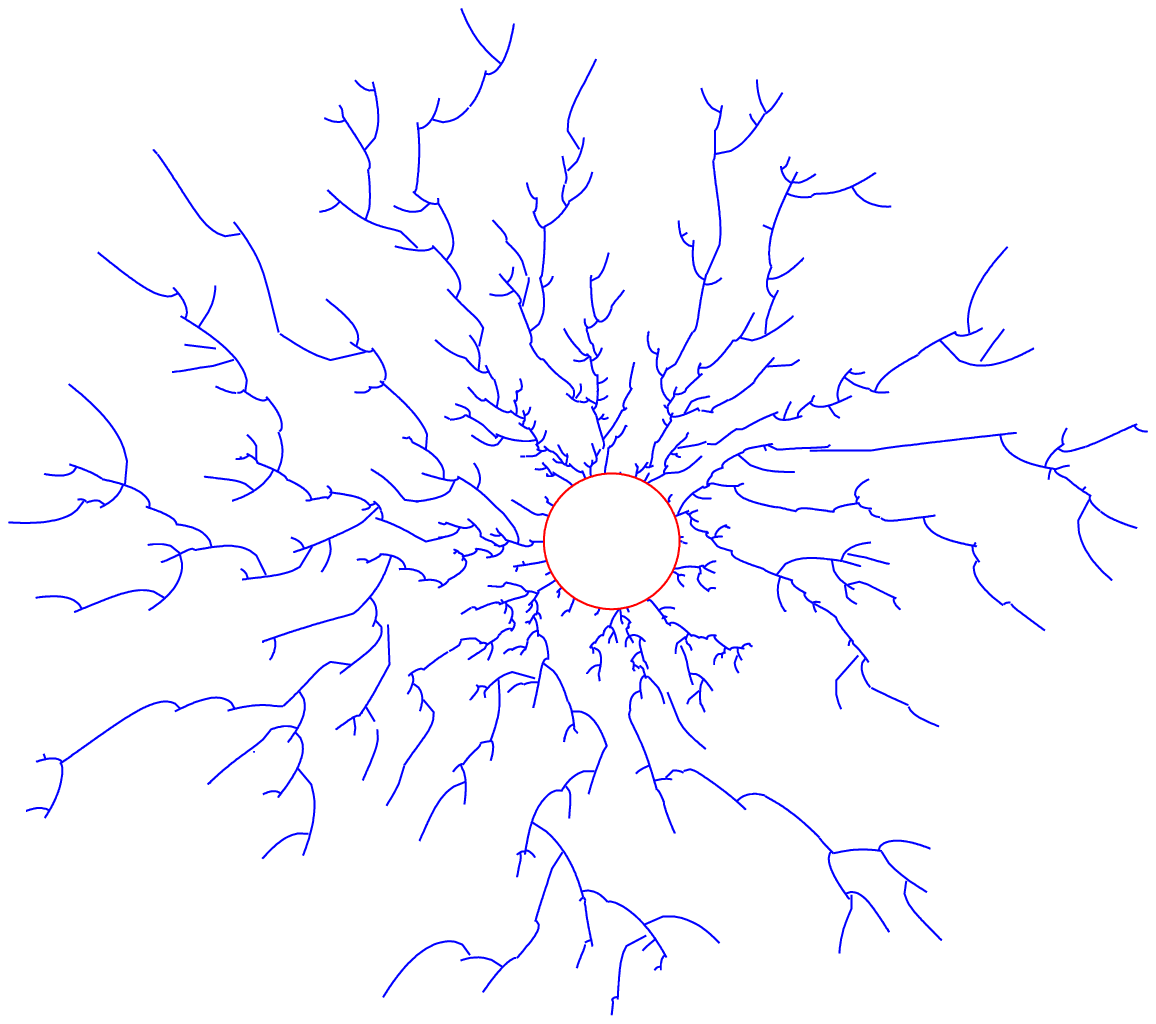}
\caption{Left: compound Poisson Loewner evolution with $r=0$, 
$\lambda=50$ at $T=2$. Right: compound Poisson Loewner evolution with 
$r=0, \lambda=300$ at $T=2$.} 
\end{figure}
We can now make our previous remarks more precise. 
We let $d_H$ denote the Hausdorff metric on compact subsets of $\CC$. 
\begin{prop} \label{det-hulls}
Fix $T,\, \lambda > 0$. Let $K_T^{r}$ be the hull
obtained by evaluating the solution to (\ref{partialloewner1}) with (\ref{xidef}) as driving function at $t=T$. 
Then, as $r \to 1$, the law of $K_T^r$ converges weakly with respect to $d_H$ to a point 
mass at the set $\overline{\DD} \cup [1, 1+\delta]$, where 
\begin{equation*}
\delta = 2e^{2T}(1+\sqrt{1-e^{-2T}})-2.
\end{equation*} 
\end{prop}
\begin{proof}
The discussion above shows that the 
driving process $\xi_{r}$ converges in law to the function identically $1$, 
when $r \to 1$. Hence we may couple the sequence $\xi_r$ so that $\xi_r \to 1$ 
a.s., see \cite[Page 85]{Durrett}.

Note that the capacity of $K_T^{r}$ is $e^T$ for all $r$. Put $\Gamma_T^r= K_T^{r} \setminus \overline{\DD}$. By Lemma 
\ref{cone} below, $\Gamma_T^{r}$ is contained in an arbitrarily narrow cone 
with apex at $0$ and with the point $1$ on its bisector. 
By the distortion estimates, we have
\begin{equation*}
e^T \le \diam(K_T^r) \le 4e^T.
\end{equation*}
We conclude that the Hausdorff limit of 
$\Gamma_T^{r}$ is $(1,1+\delta]$ for some $\delta$. Note that
the limiting function $f_T:\Delta \to \Delta \setminus [1,1+\delta]$ has 
expansion $e^Tz+\cdots$ at infinity. This together with the 
explicit formula for $f_T$ yields 
the stated expression for $\delta$. 
\end{proof}
The following is the radial version of Lemma 2.1 from \cite{Rohde_Stable}.
\begin{lem}
\label{cone}
Let $\xi(t)$ be the driving function in (\ref{partialloewner1}). Let
 $\theta_1 < \theta_2$ where $\theta_{j} \in [-\pi, \pi), \, j=1,2,$ and suppose that $\arg \xi(t) \in [\theta_1, \theta_2]$ for $0 \le t \le T$. Then
\begin{equation*}
K_T\setminus \overline{\mathbb{D}}
 \subset \{re^{i \theta}: \theta_1 \le \theta \le \theta_2,\, r \geq 0\}.
\end{equation*}
\end{lem}
We shall now drop the assumption that $\lambda$ be kept fixed and consider the 
case when $\lambda=\lambda(r)=1/(1-r)$. We write
$\xi_r=\exp(i Y_r(t))$ for the corresponding unimodular driving process. In 
this case the Fourier coefficients are
\begin{equation}
a_{Y_r(t)}(n)=\exp \left[t \fr{ r^{|n|} - 1}{1-r} 
\right],\quad n \neq 0.
\label{Ycoeffs_2}
\end{equation} 
Letting $r \to 1$ in (\ref{Ycoeffs_2}), we observe that 
$(r^{|n|} - 1)/(1-r) \to -|n|$.
Hence
\begin{equation*}
\lim_{r \to 1}a_{Y_r(t)}(n)=\exp \left[-t|n|\right],\quad n \neq 0,
\end{equation*} 
and these are the Fourier coefficients of a Cauchy process on $\TT$. By 
appealing to Corollary \ref{continuity} and Lemma \ref{levy-conv} we 
have proved
\begin{prop}\label{cauchy}
Fix $T >0$. For each $r< 1$ let $f_T^r$ denote the solution to 
(\ref{partialloewner1}) at $t=T$ with $\xi_r=\xi_{r,1/(1-r)}$ as driving process.
Let $f$ denote the solution to (\ref{partialloewner1}) at $t=T$ with a unimodular 
Cauchy process as driving process.
Then, as $r \to 1$, the law of $e^{-T}f_T^{r}$ converges weakly in 
$\Pi(\Sigma')$ to that of $e^{-T}f_T$.
\end{prop} 
A systematic study of Loewner evolutions in a half-plane driven by the 
Cauchy process has recently been carried out by Chen and Rohde (see 
\cite{Rohde_Stable}).

\subsection{Related models}
We have also studied a related model. Instead of keeping 
the parameter $r$ fixed in the definition of the law of $X^r$, we
let $r$ be a random variable on $[0,1]$ with density
\begin{equation*}
f_r^{\beta}(x)=
\beta x^{\beta-1} \quad \rm{for} \quad \beta > 0. 
\end{equation*}
A computation then shows that the density of the new random variable 
$X_{\beta}$ is given by
\begin{equation} 
f_{X_\beta}(\theta)=\frac{1}{2\pi}+\frac{\beta}{\pi}\sum_{n=1}^{\infty} 
\frac{\cos n\theta}{n+\beta},\quad \theta \in [-\pi,\pi].
\label{densfunct} 
\end{equation}
We then define the compound Poisson process $Y_{\beta, \lambda}$ as 
before, and use $\xi_{\beta,\lambda}=\exp(i Y_{\beta,\lambda})$ as 
driving function. The Fourier coefficients of the process
are
\begin{equation*}a_{Y(t)}(n)=\exp\left\{-\lambda t\left(\frac{n}{\beta+n}\right)\right\}.\end{equation*}
The qualitative behavior of this model is similar, and all 
our results go through with minor modifications: $\beta$ is a localizing 
parameter, and we obtain the Cauchy process 
in the limit $\beta=\lambda$, $\beta \rightarrow \iy$.
\begin{figure}
\centering
\includegraphics[height=0.24 \textheight, width=0.48 \textwidth]{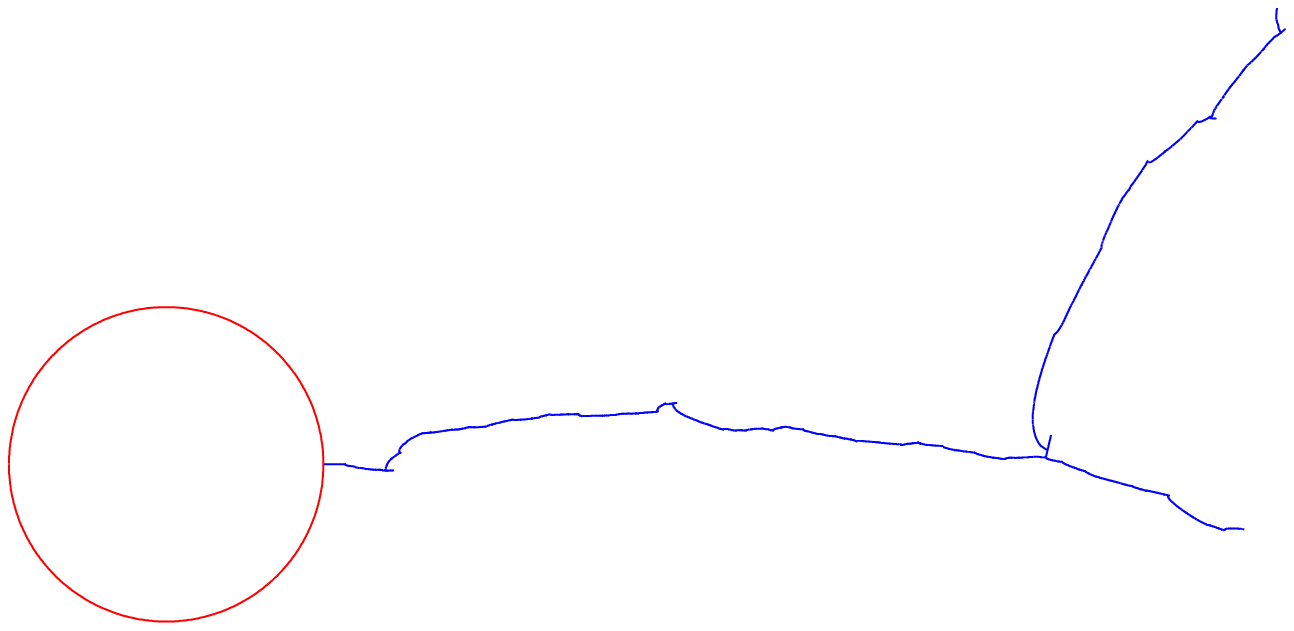}
\includegraphics[height=0.24 \textheight, width=0.48 \textwidth]{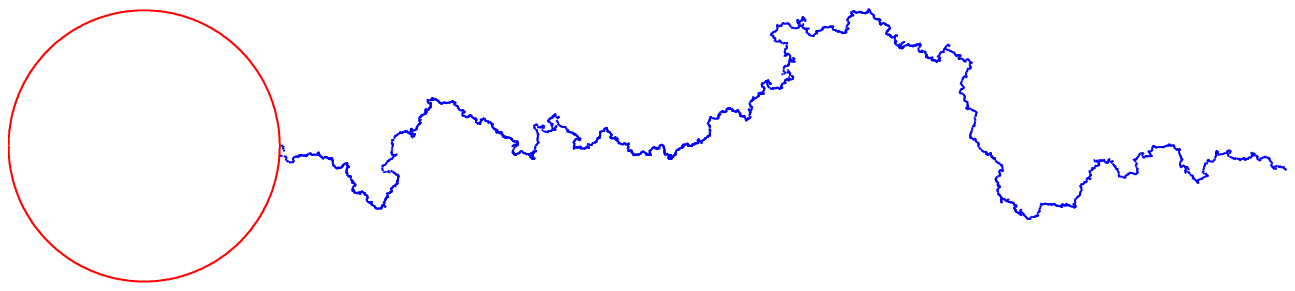}
\caption{Left: compound Poisson Loewner evolution with 
$r=\beta= \lambda=1000$ at $T=1$. Right: compound Poisson Loewner evolution with 
$\gamma=8\pi^2/ \lambda$, 
$\lambda=4000$ at $T=1$.} 
\end{figure}
We could of course exchange the Poisson kernel in the definition of 
the density of the jumps of the compound Poisson process in favor 
of some other density. For instance, we could choose $X_{\gamma}$ 
distributed according to the heat kernel with parameter $\gamma>0$,
\begin{equation*}k_{\gamma}(\theta)=\sum_{n=-\infty}^{\infty}e^{-\gamma n^2}
e^{i n\theta}.\end{equation*}
In this case, there is no interpretation of the growth model in 
terms of harmonic measure, except if we let $\gamma \rightarrow \iy$ since 
we then recover our previous model with $r=0$. However, we may note the 
following. The Fourier coefficients of the resulting process 
$\xi=\exp(i Y_{\gamma,\lambda})$ are
\begin{equation*}a_{Y(t)}(n)=\exp\{\lambda t[\exp(-n^2 \gamma)-1]\}.\end{equation*} 
Taking $\lambda=c/\gamma$ and letting $\gamma \to 0$, we see that 
the driving process converges to Brownian motion on the unit circle. Hence, 
by Lemma \ref{levy-conv},
the laws of the conformal mappings converge weakly to the 
law of the $\textrm{SLE}(\kappa)$ mapping if we 
choose $c=4\pi^2\kappa$. This result is analogous to 
\cite[Proposition 2]{Bauer}, except that Bauer considers a different space of 
conformal mappings.
\subsection{Hausdorff dimension of the limiting rescaled hulls}
We again turn our attention to the model related to the Poisson kernel, where
the driving function is given by (\ref{xidef}).
In this section, the dependence of the hulls on the parameter $\lambda$
will be of importance and, as before, we indicate this by writing 
$K_t^{r,\lambda}$ for the hulls of the evolution. The limit 
\begin{equation*}
\tilde{K}^{r,\lambda}_{\infty}=
\lim_{t \rightarrow \iy} e^{-t}K_t^{r,\lambda}
\end{equation*}
exists in the sense of Theorem \ref{existence}, and we denote the 
corresponding limit measure on $\Sigma'$ by $\PP^{r,\lambda}_{\iy}$. 

The following theorem shows that the Hausdorff 
dimension of the limit hulls is trivial, just as in the $\rm{HL}(0)$ case.
The proof follows along the main lines of Section 5 in \cite{Rohde_Zins}.
The Hausdorff and upper box-counting dimensions of a set will be denoted by 
$\dim_H$ and $\overline{\dim}_B$, respectively (see \cite{Mattila} for 
definitions).
\begin{thm}
Let $r=0$, and fix $\lambda>4$. Then 
\begin{equation*}
\dim_H(\partial \tilde{K}^{0,\lambda}_{\iy})=1,
\label{Hdim}
\end{equation*}
almost surely.
\end{thm}
\begin{proof}
We write $K_t=K_t^{0,\lambda}$.
Let $\{\tau_k\}_k$ denote the waiting times of the Poisson 
process $N$. We note that the $\tau_k$:s are independent and exponentially
distributed with density $\lambda e^{-\lambda t}$. Next, we write
\begin{equation*}
T_n=\sum_{k=1}^{n}\tau_k, \quad n \geq 1.
\end{equation*}
We find it convenient to consider the hulls 
$K_n=K_{T_n}$ and their rescalings; clearly, 
$\tilde{K}_n \rightarrow \tilde{K}_{\iy}$ as $n \rightarrow \iy$. Also, we let
$\mathcal{F}_t=\sigma(\xi(s), \,0\leq s \leq t)$ be the natural filtration and
write $\mathcal{F}_n=\mathcal{F}_{\tau_n}$ for the stopped sigma algebras.

The hulls $K_n$ may be written inductively as the union of the unit disk and
$n$ curvilinear arcs $l_k=K_k\setminus K_{k-1}$:
\begin{equation*}K_n=\overline{\DD} \cup \bigcup_{k=1}^{n}l_k,\end{equation*}
We denote the length of each arc $l_k$ by $\ell_k$.
The length of $\partial K_n$ is then given by
\begin{equation*}L_n=2\pi+\sum_{k=1}^n \ell_k.\end{equation*} 
Similarily, the length of $\partial \tilde{K}_n$ is
\begin{equation*}\tilde{L}_n=e^{-T_n}\left(2\pi+ \sum_{k=1}^n\ell_k\right).\end{equation*}
Our strategy is to show that $\tilde{L}_{\iy}$ is finite almost surely, and this will imply the desired conclusion.
To do this we prove a uniform bound on the expected lengths $\EE[\tilde{L}_n]$ of 
the hulls $\tilde{K}_n$. 

We note that
\begin{align}
\EE[\tilde{L}_n]&=2\pi \EE[e^{-T_n}]+\sum_{k=1}^n \EE[e^{-T_n}\ell_k]\nonumber
\\
&\leq 2\pi (\EE[e^{-2T_n}])^{\frac{1}{2}}+\sum_{k=1}^n
(\EE[e^{-2T_n}\ell_k^2])^{\frac{1}{2}},
\label{CSsum}
\end{align}
and thus, it suffices to estimate $\EE[e^{-2T_n}\ell_k^2]$ for 
$k=1,\ldots,n$.

Let $f_{k-1}:\Delta \rightarrow \mathbb{C}\setminus K_{k-1}$ 
denote the solution at time $T_{k-1}$ to (\ref{partialloewner})
with the compound Poisson driving function $\xi$. The 
length of the $k^{th}$ added arc is given by the random variable
\begin{equation*}\ell_k=\int_{1}^{1+\delta_k}|f'_{k-1}(re^{i\theta})|dr,\end{equation*}
where $\theta$ is chosen uniformly. By the Cauchy-Schwarz inequality,
\begin{equation*}\ell_k^2\leq \delta_k\int_1^{1+\delta_k}|f'_{k-1}(r e^{i\theta})|^2dr.\end{equation*}
Since $\theta$ is uniformly distributed, an estimate of 
$\EE[e^{-2T_n}\ell_k^2]$ requires that we estimate the integral
\begin{equation}
\frac{1}{2\pi}\int_{-\pi}^{\pi}\int_1^{1+\delta_k}
|f'_{k-1}(re^{i\theta})|^2dr d\theta.
\label{rthetaint}
\end{equation}
We note that $f_{k-1}(z)=e^{T_{k-1}}z+\cdots $, and apply 
Lemma \ref{intestim} below to obtain
\begin{equation*}
\frac{1}{2\pi}\int_{-\pi}^{\pi}\int_1^{1+\delta_k}
|f'_{k-1}(re^{i\theta})|^2dr d\theta\leq C e^{2T_{k-1}} (1+\delta_k)^2,
\end{equation*}
where $C<\iy$ is a constant that does not depend on $k$.
Hence, 
\begin{equation}
\EE[e^{-2\tau_k}\ell_k^2|\mathcal{F}_{k-1}]\leq
 C e^{2T_{k-1}}\EE[e^{-2\tau_k}\delta_k\, (1+\delta_k)^2]
\label{condestim}
\end{equation}
holds almost surely.
Using (\ref{condestim}) and properties of conditional expectation
and independence, we obtain
\begin{align}
\EE[e^{-2T_n}\ell_k^2]&
=\EE[e^{-2(\tau_{k+1}+\cdots+\tau_n)}]\,\EE[e^{-2T_{k}}\ell_k^2]\\
&=(\EE[e^{-2\tau}])^{n-k}\,
\EE\left[ \EE[e^{-2T_{k}}\ell_k^2|\mathcal{F}_{k-1}]\right] \nonumber\\
&=(\EE[e^{-2\tau}])^{n-k}\,
\EE\left[e^{-2T_{k-1}} \EE[e^{-2\tau_k}\ell_k^2|\mathcal{F}_{k-1}]\right]
 \nonumber\\
&\leq
C\cdot (\EE[e^{-2\tau}])^{n-k}\, \EE[e^{-2\tau_k}\delta_k \, (1+\delta_k)^2] \nonumber
\end{align}
A computation shows that $\EE[e^{-2\tau_k}]=\lambda/(\lambda+2)$ and
in view of (\ref{deltadist}), we have $\delta_k<4 e^{2\tau_k}$.  
Hence $\EE[e^{-2\tau_k}\delta_k \,(1+\delta_k)^2]< \iy$ because 
of our assumption $\lambda>4$. This 
leads to the estimate
\begin{equation*}\EE[e^{-2T_n}\ell_k^2] \leq C \cdot
\left(\frac{\lambda}{\lambda+2}\right)^{n-k},\end{equation*}
where $C=C(\lambda)$, which in turn implies that
\begin{equation}
\EE[\tilde{L}_n]\leq C\sum_{k=0}^n
\left(\frac{\lambda}{\lambda+2}\right)^{\frac{n-k}{2}}\\
\leq C
\sum_{k=0}^{\iy}\left(\frac{\lambda}{\lambda+2}\right)^{\frac{k}{2}}. 
\end{equation}
Since $\lambda/(\lambda+2)<1$, the geometric series on the right-hand side
converges, and we conclude that $\EE[\tilde{L}_n]$ is bounded by a
constant that depends on $\lambda$ but not on $n$, which in turn 
implies that 
\begin{gather}
\label{L-bound}
\EE[\tilde{L}_{\iy}] \le C
\end{gather}
for some constant $C=C(\lambda) < \iy$.

Let $N(\delta, \partial \tilde{K}_{\iy})$ denote the minimal number of 
disks of radius $\delta$ needed to cover $\partial \tilde{K}_{\iy}$. There exists a constant $C<\iy$ independent of $\delta$ such
that, almost surely, 
\begin{equation*}
N(\delta, \partial \tilde{K}_{\iy}) \le C \tilde{L}_{\iy}/\delta < \iy.
\end{equation*} 
Hence, we have $\dim_H(\partial \tilde{K}_{\iy})\leq 
\overline{\dim}_B(\partial \tilde{K}_{\iy})\leq 1$, and the proof is 
complete.
\end{proof}
\begin{cor}
Fix $r<1$ and $\lambda>4$. Then 
\begin{equation*}
\dim_H(\partial \tilde{K}^{r,\lambda}_{\iy})=1,
\end{equation*}
almost surely.
\end{cor}
\begin{proof}
It is enough to note that the density of $X^r$ is square integrable when 
$r \neq 0$, and to use this fact together with the Cauchy-Schwarz 
inequality in (\ref{CSsum}). From there, we proceed as in the previous proof.
\end{proof}
\begin{lem}\label{intestim}
For $f \in \Sigma'$ and $\delta>0$, we have
\begin{equation*}
\frac{1}{2\pi}\int_{-\pi}^{\pi}\int_1^{1+\delta}
|f'(re^{i\theta})|^2dr d\theta
\leq C (1+\delta)^2,
\end{equation*}
where $C<\infty$ is a universal constant.
\end{lem}
\begin{proof}
We use the standard technique of inversion to obtain a function 
$F(\zeta)=1/f(1/\zeta)$, $\zeta \in \DD$,
belonging to the class $\mathcal{S}$ of univalent 
functions in the unit disk. A computation shows that 
$f'(z)=F'(1/z)/(z^2(F(1/z))^2)$.
We plug this into the integral and change variables 
to obtain
\begin{equation*}
\frac{1}{2\pi}\int_{-\pi}^{\pi}\int_1^{1+\delta}
|f'(re^{i\theta})|^2dr d\theta
=\frac{1}{2\pi}\int_{-\pi}^{\pi}\int_{(1+\delta)^{-1}}^1
\frac{|F'(\rho e^{i\eta})|^2}{|F(\rho e^{i \eta})|^4}\rho^2d\rho d\eta.
\end{equation*}
The integral on the right-hand side is bounded by
\begin{equation*}
\frac{1}{2}\int_{\mathcal{A}((1+\delta)^{-1},1)}
\frac{|F'(\rho e^{i\eta})|^2}{|F(\rho e^{i \eta})|^4} dA(w),
\end{equation*}
where $dA=\rho d\rho d\eta/ \pi$ denotes area measure and 
\begin{equation*}
\mathcal{A}((1+\delta)^{-1},1))=\{w \in \DD:(1+\delta)^{-1}<|w|<1\}.
\end{equation*}  
Another change of variables yields
\begin{equation*}
\frac{1}{2}\int_{\mathcal{A}((1+\delta)^{-1},1)}
\frac{|F'(\rho e^{i\eta})|^2}{|F(\rho e^{i \eta})|^4} dA(w)
=\frac{1}{2}\int_{F\left(\mathcal{A}((1+\delta)^{-1},1)\right)}|z|^{-4}dA(z),
\end{equation*}
and by applying growth estimates for $\mathcal{S}$, 
\begin{equation*}
|F(z)|\geq \frac{(1+\delta)^{-1}}{(1+(1+\delta)^{-1})^2}, 
\quad z \in \mathcal{A}((1+\delta)^{-1},1),
\end{equation*}
we find that the last integral is bounded by $C\cdot(2+\delta)^4/(1+\delta)^2<
C'\cdot (1+\delta)^2$.
\end{proof}
\begin{rem}
The hypothesis $\lambda>4$ is only used to obtain an easy estimate of 
$\EE[e^{-2\tau_k}\delta_k \,(1+\delta_k)^2]$, and is probably not necessary for 
the conclusion since small $\lambda$ correspond to long arcs being added.  
\end{rem}
\begin{rem}
We note that our result on Hausdorff dimension of the rescaled hulls is also 
valid for the other compound Poisson process-driven evolutions discussed in 
the previous subsection. 
\end{rem}
\subsection*{Acknowledgements}
We thank H\aa kan Hedenmalm for suggesting the model of 
section 4 and for numerous useful comments and suggestions. 
We also thank Michael Benedicks and Boualem Djehiche for interesting 
discussions on the topics in this paper. We are grateful to 
the anonymous referee whose detailed comments helped us improve the 
quality of our paper.
\bibliography{js}
\end{document}